\documentclass[12pt]{article}
\usepackage{mathrsfs}
\usepackage{amssymb}
\usepackage{amsfonts}
\usepackage{color,xcolor}
\usepackage{graphicx}
\usepackage{manfnt}
\usepackage[dvips,  
pdfstartview=FitH, CJKbookmarks=true,
bookmarksnumbered=true, bookmarksopen=true,  colorlinks, 
pdfborder=001,   
linkcolor=red,
anchorcolor=red,
citecolor=blue]{hyperref}
\newtheorem{theorem}{Theorem}[section]
\newtheorem{corollary}{Corollary}[section]
\newtheorem{lemma}{Lemma}[section]

\newtheorem{definition}{Definition}[section]
\newtheorem{remark}{Remark}[section]

\def\[{{\Big[}}\def\]{{\Big]}}\def\<{{\langle}}\def\>{{\rangle}}\def\({{\Big(}}
\def\){{\Big)}}

\def\min{{\mathord{{\rm min}}}}
\def\={&\!\!=\!\!&}
\def\cB{{\mathcal B}}\def\cC{{\mathcal C}}
\def\cF{{\mathcal F}}

\def\cM{{\mathcal M}}

\def\mE{{\mathbb E}}

\def\mN{{\mathbb N}}\def\mP{{\mathbb P}}
\def\mR{{\mathbb R}}

\def\geq{\geqslant}\def\leq{\leqslant}

\begin{document}
\title{\bf Schauder estimates for stochastic transport-diffusion equations with L\'{e}vy processes}
\author{Jinlong Wei$^a$, Jinqiao Duan$^b$  and Guangying Lv$^c$}
\date{{\it $^a$School of Statistics and Mathematics, Zhongnan University of}\\
{ \it Economics and Law, Wuhan, Hubei 430073, China}\\
{ \tt  weijinlong@zuel.edu.cn}
\\
{\it $^b$Department of Applied Mathematics}\\
{\it Illinois Institute of Technology, Chicago, IL 60616, USA}\\
{ \tt duan@iit.edu}\\
{\it $^c$Institute of Contemporary Mathematics, Henan University}\\
{\it Kaifeng, Henan 475001, China}\\
{\tt gylvmaths@henu.edu.cn}}

\maketitle{}
\noindent{\hrulefill}
\vskip1mm\noindent
{\bf Abstract} We consider a transport-diffusion equation with
L\'{e}vy noises and H\"{o}lder continuous coefficients. By using the heat kernel estimates, 
we derive the Schauder estimates for the mild solutions. Moreover, 
when the transport term vanishes and $p=2$, we show that the H\"{o}lder index in space variable is optimal.



 \vskip0mm\noindent{\hrulefill}
\section{Introduction}\label{sec1}\setcounter{equation}{0}
Let $(\Omega,\cF,\{\cF_t\}_{t\geq0},\mP)$ be a filtered complete probability space
with the right continuous filtration $\cF_t$.
 Denote $\{W_t\}_{t\geq0}$   a  scalar Wiener process on $(\Omega,\cF,\{\cF_t\}_{t\geq0},\mP)$. Let $E$ be a ball $B_c(0)-\{0\}$,  of radius $c$ without the center. Moreover,  $\tilde{N}$ be a time homogeneous compensated Poisson random measure defined on $(\Omega,\cF,\{\cF_t\}_{t\geq0},\mP)$ (defined in Definition \ref{def2.2}), which is independent of $\{W_t\}_{t\geq0}$ and has an intensity measure $\nu\times \lambda$ on $E\times\mR_+$.

In the present paper, we are concerned with the existence, uniqueness and regularity of  the mild solution for the following stochastic transport-diffusion equation:
\begin{eqnarray}\label{1.1}
&&du(t,x)-b(t,x)\cdot\nabla u(t,x)dt-\frac{1}{2}\Delta u(t,x)dt
\cr&&=h(t,x)dt+f(t,x)dW_t+\int\limits_Eg(t,x,v)\tilde{N}(dt,dv), \ \ t>0, \ x\in \mR^d.
\end{eqnarray}

\medskip

 When the L\'evy noise part absent ($g=0$),  this stochastic partial differential equation (SPDE)  (\ref{1.1}) has been studied widely. When $g=h=0$, $b=0$ and the initial datum vanishes, (\ref{1.1}) becomes
\begin{eqnarray}\label{1.2}
du(t,x)-\frac{1}{2}\Delta u(t,x)dt=f(t,x)dW_t, \ \ t>0, \ x\in \mR^d, \quad u(t,x)|_{t=0}=0.
\end{eqnarray}
For example, Krylov  \cite{Kry99} obtained  the following estimate for the solution of the Cauchy problem (\ref{1.2})  for $p\geq 2$,
\begin{eqnarray}\label{1.3}
\mE\|\nabla u\|_{L^p((0,T)\times\mR^d)}^p\leq C(d,p)
\mE\|f\|_{L^p((0,T)\times\mathbb{R}^d)}^p,
\end{eqnarray}
using a variant of   the   Littlewood-Paley inequality.
This result was extended by Neerven, Veraar and Weis \cite{NVW}  to the case when the Laplace operator $\Delta$ is replaced by a   linear operator  $A$ which admits a bounded $H^\infty$-calculus of angle less than $\pi/2$. For research       in the $L^p$-theory for linear SPDEs, see \cite{Kim08,Kry96,Kry00,Kry08} and  in the $L^p$-theory for nonlinear SPDEs, refer to  \cite{DHV,DMH,KK16,Zha}.
For $p=\infty$, an analogue and interesting estimate of (\ref{1.3}) for (\ref{1.2}) was also derived by Denis, Matoussi and Stoica \cite{DMS}. By using Moser's iteration scheme developed by Aronson and Serrin,   they derived  a space-time $L^\infty$ estimates for certain nonlinear SPDEs. Moreover, after introducing a notion of stochastic BMO spaces,  Kim \cite{Kim15} obtained a BMO estimate for $\nabla u$, which is controlled by $\|f\|_{L^\infty}$.    For more details in this topic, one also sees \cite{Kuk,LGWW}.

There exist also some Schauder estimates for solutions of  (\ref{1.1}) when L\'evy noise is absent ($g=0$).  When $f(t,\cdot)$ belongs to $L^p$ with sufficiently large $p$ (or $p=\infty$) and $\mR^d$ is replaced by a bounded domain $Q$ (with smooth boundary), the time and space $\cC^\alpha$ estimates have been discussed by Kuksin, Nadirashvili and Piatnitski \cite{KNP}.   This result was further strengthened by \cite{Kim04},  for  a general H\"{o}lder estimates for generalized solutions with $L_q(L_p)$ coefficients.  Later, Du and Liu \cite{DL} extended the result on bounded domains to $\mR^d$ and built the $C^{2+\alpha}$-theory. When $f$ and $h$ are dependent on $u$ (nonlinear SPDE case), the $\cC^\alpha$ estimates were also derived by Hsu, Wang and Wang \cite{HWW}. They use a stochastic De Giorgi iteration technique and  proved that the solution is almost surely $\cC^\alpha$   in both space and time. When $u$ takes values in a Hilbert space, some regularity results are available  \cite{DMH, Mik}.

\medskip

When the L\'evy noise part is present ($g \neq 0$),  Kotelenez \cite{Kot84}, and  Albeverioa, Wu and Zhang \cite{AWZ} studied the $L^2$-theory for  the SPDE (\ref{1.1}). Moreover,   an   $L^p$ theory is founded by Marinelli, Pr\'{e}v\^{o}t and R\"{o}ckner \cite{MPR}.

However, as far as we know, there have been very few papers   dealing with the Schauder estimates for  (\ref{1.1}). In this  present paper, we will fill this gap and derive the Schauder estimates for  the mild solutions.

\medskip

This paper is organized as follows. After introducing some notions and stating the main result in Section 2, we present several useful lemmas in Section 3. In Section 4, we prove the main result.    Finally, we conclude  with some remarks   on the regularity of mild solutions to problem (\ref{1.1}) in Section \ref{sec5}.

\vskip2mm\noindent
\textbf{Notations} Denote $B_r(x):=\{y\in\mathbb{R}^d:|x-y|<r\}$ by the ball centered at $x$ with radius $r$.  $a\wedge b=\min\{a,b\}$, $a\vee b=\max\{a,b\}$. $\mR_+=\{r\in \mR, \ r\geq 0\}$. The letter $C$ will mean a positive constant, whose values may change in different places. The  Lebesgue measure is denoted by  $\lambda$,  or  by $dt$ if there is no confusion. $\mN$ is  the  set of natural numbers. Let $\mN_0:=\mN\cup \{0\}$ and $\overline{\mN}:=\mN_0\cup \{\infty\}$. $\cB(E)$ is the Borel $\sigma$-algebra on $E$.   By $M_+(E)$ we denote the family of all $\sigma$-finite positive measures on $E$ , by $\cM_+(E)$ we denote the $\sigma$-field on $M_+(E)$ generated by functions $i_B: M_+(E)\ni\mu\rightarrow \mu(B)\in\mR_+, B\in\cB(E)$.

\section{Main result}\label{sec2}
\setcounter{equation}{0}
Let $E$ and $(\Omega,\cF,\{\cF_t\}_{t\geq0},\mP)$ be given in the previous section. We first recall the notion of Poisson random measure.

\begin{definition} \label{def2.1}  A time homogeneous Poisson random measure $N$ on $(E,\cB(E))$ over the filtered probability space $(\Omega,\cF,\{\cF_t\}_{t\geq0},\mP)$ with an intensity measure $\nu\times \lambda$, is a measurable function $N: (\Omega, \cF)\rightarrow(M_+(E\times\mR_+), \cM_+(E\times\mR_+))$, such that

(i) for each $B\times I \in \cB(E)\times  \cB(\mR_+)$, if $\nu(B)<\infty$, $N(B\times I)$ is a Poisson random variable with parameter $\nu(B)\lambda(I)$;

(ii) $N$ is independently scattered, i.e. if the sets $E_j\times I_j\in \cB(E)\times  \cB(\mR_+), \ j=1,_{\cdots},n$ are pairwise disjoint, then the random variables $N(B_j\times I_j), \  j=1,_{\cdots},n$ are mutually independent; and

(iii) for each $U\in \cB(E)$, the $\overline{\mN}$-valued process $\{N((0,t],U)\}_{t>0}$ is $\{\cF_t\}_{t\geq0}$-adapted and its increments are independent of the past.
\end{definition}

\begin{remark} \label{rem2.1}
In this definition, $\nu$ is called a L\'{e}vy measure and it satisfies the following condition
$$
\int_E1\wedge |v|^2\nu(dv)<\infty.
$$
\end{remark}

\begin{definition} \label{def2.2} Let $N$ be a homogeneous Poisson random measure on $(E,\cB(E))$ over the probability space $(\Omega,\cF,\{\cF_t\}_{t\geq0},\mP)$. The $\mR$-valued process $\{\tilde{N}((0,t],A)\}_{t>0}$ defined by
$$
\tilde{N}((0,t],A)=N((0,t],A)-\nu(A)t, \quad t>0, \ A\in \cB(E),
$$
is called a compensator Poisson random measure. And now $\{\tilde{N}((0,t],A)\}_{t>0}$ is a martingale on $(\Omega,\cF,\{\cF_t\}_{t\geq0},\mP)$.
\end{definition}

In this paper, our focus will be on Schauder estimates of mild solutions for (\ref{1.1}). To formulate the Cauchy problem, we assume that the initial value vanishes. Here the mild solution is defined as follows:

\begin{definition} \label{def2.3} Let $u$ be a $\cB(\mR_+)\times\cB(\mR^d)\times\cF$ measurable function. We call that $u$ is a mild solution of (\ref{1.1}), with initial data vanishes, if the following properties hold:

(1) $u$ is $\cF_t$-adapted;

(2) $\{u(t,x,\cdot)\}_{t\geq 0}$ as a family of $L^2(\Omega,\cF,\mP)$-valued random variables is right continuous and has left limits in the variable $t\in [0,\infty)$, namely,
\begin{eqnarray}\label{2.1}
u(t-,x,\cdot)=L^2(\Omega)-\lim_{s\uparrow t}u(s,x,\cdot), \ t\in [0,\infty);
\end{eqnarray}

(3) $u \in L^\infty_{loc}([0,\infty);W^{1,\infty}(\mR^d;L^2(\Omega)))$;

(4) for every $t>0$, the following equation
\begin{eqnarray}\label{2.2}
u(t,x)&=&\int\limits_0^t P_{t-r} (b(r,\cdot)\cdot\nabla u(r,\cdot))(x)dr+\int\limits_0^t P_{t-r}h(r,\cdot)(x)dr+\int\limits_0^t P_{t-r}f(r,\cdot)(x)dW_r\cr&&+\int\limits_{(0,t]} \int\limits_EP_{t-r}g(r,\cdot,v)(x)\tilde{N}(dr,dv),
\end{eqnarray}
holds almost surely, where the stochastic integral in (\ref{2.2}) is interpreted in It\^{o}'s and $P_t$ denotes the forward heat semigroup, i.e.
\begin{eqnarray}\label{2.3}
P_t\varphi(x)=\frac{1}{(2\pi t)^{d/2}}\int\limits_{\mR^d}e^{-\frac{|x-y|^2}{2t}}\varphi(y)dy, \quad \varphi\in L^\infty(\mR^d).
\end{eqnarray}
\end{definition}

\begin{remark} \label{rem2.2} The definition here is inspired by Marinelli, Pr\'{e}v\^{o}t and R\"{o}ckner \cite[Definition 2.1]{MPR} and the definition in \cite{AWZ}.
\end{remark}

Before stating our main result, we recall  some notations for function spaces. For   $T>0$, $\alpha>0$ and  $p\geq 2$,     define $L^\infty([0,T];\cC^{\alpha}_{b}(\mR^d))$
to be the set of all $\cC^{\alpha}_{b}(\mR^d)$-valued essentially bounded functions $u$ such that
\begin{eqnarray*}
\|u\|_{T,\infty,\alpha}:=\mathop {\mbox{esssup}}_{0\leq t\leq T}\max_{x\in\mR^d}|u(t,x)|+\mathop {\mbox{esssup}}_{0\leq t\leq T}\sup_{x,y\in\mR^d, x\neq y}\frac{|u(t,x)-u(t,y)|}{|x-y|^\alpha}<\infty.
\end{eqnarray*}
When $\alpha=0$, $\|u\|_{T,\infty,0}$ is written by $\|u\|_{T,\infty}$ for short and $\|u\|_{\infty}:=\max_{x\in\mR^d}|u(x)|$. Correspondingly, $L^\infty([0,T];\cC^{1,\alpha}_{b}(\mR^d))$ is the set of all functions in $L^\infty([0,T];\cC^{\alpha}_{b}(\mR^d))$, such that
\begin{eqnarray*}
\|u\|_{T,\infty,1+\alpha}&:=&\mathop {\mbox{esssup}}_{0\leq t\leq T}\max_{x\in\mR^d}|u(t,x)|+\mathop {\mbox{esssup}}_{0\leq t\leq T}\max_{x\in\mR^d}|Du(t,x)|\cr&&+\mathop {\mbox{esssup}}_{0\leq t\leq T}\sup_{x,y\in\mR^d, x\neq y}\frac{|Du(t,x)-Du(t,y)|}{|x-y|^\alpha}<\infty.
\end{eqnarray*}
Similarly,   define the spaces $L^\infty([0,T];\cC^{\alpha}_{b}(\mR^d;L^p(\Omega))), L^\infty([0,T]; L^p(E,\nu;\cC^{\alpha}_{b}(\mR^d)))$. For $h\in L^\infty([0,T];\cC^{\alpha}_{b}(\mR^d;L^p(\Omega)))$ and $g\in L^\infty([0,T];L^p(E,\nu;\cC^{\alpha}_{b}(\mR^d)))$, the norms are given by
\begin{eqnarray*}
\|h\|_{T,\infty,\alpha,p}:=\mathop {\mbox{esssup}}_{0\leq t\leq T}\max_{x\in\mR^d}\|h(t,x)\|_{L^p(\Omega)}+\mathop {\mbox{esssup}}_{0\leq t\leq T}\sup_{x,y\in\mR^d, x\neq y}\frac{\|h(t,x)-h(t,y)\|_{L^p(\Omega)}}{|x-y|^\alpha}<\infty
\end{eqnarray*}
and
\begin{eqnarray*}
\|g\|_{T,\infty,p,E,\alpha}&:=&
\mathop {\mbox{esssup}}_{0\leq t\leq T}\|\max_{x\in\mR^d}|g(t,x,\cdot)|\|_{L^p(E,\nu)}
\cr&&+\mathop {\mbox{esssup}}_{0\leq t\leq T}\Big \|\sup_{x,y\in\mR^d, x\neq y}\frac{|g(t,x,\cdot)-g(t,y,\cdot)|}{|x-y|^\alpha} \Big \|_{L^p(E,\nu)}<\infty
\end{eqnarray*}
respectively.

Our main result is as follows.
\begin{theorem} \label{the2.1} Let $b,h,f$ and $g$ be measurable functions. Consider the stochastic transport-diffusion equation (\ref{1.1})   with the  zero initial data. For $\alpha>0$, $p>2$, we assume that
\begin{eqnarray}\label{2.4}
0<\alpha +\frac{2}{p}-1=:\gamma,
\end{eqnarray}
and
\begin{eqnarray}\label{2.5}
f \in L^\infty_{loc}([0,\infty);\cC^{\alpha}_{b}(\mR^d)), \ g\in L^\infty_{loc}([0,\infty);L^{p+}(E,\nu;\cC^{\alpha}_{b}(\mR^d))) \ g(t,x,\cdot) \ \mbox{vanishes near} \ 0.
\end{eqnarray}
In addition, we assume that there is a real number $0<\beta<\gamma$, such that
\begin{eqnarray}\label{2.6}
b\in L^\infty_{loc}([0,\infty);\cC^{\beta}_{b}(\mR^d;\mR^d)), \ h\in L^\infty_{loc}([0,\infty);\cC^{\beta}_{b}(\mR^d;L^p(\Omega))).
\end{eqnarray}
Then there exists  a unique mild solution $u$ for the equation  (\ref{1.1}). Moreover,
$u$ is in class of  $L^\infty_{loc}([0,\infty);\cC^{1+\gamma-}_b(\mR^d;L^p(\Omega)))$
and for every $t>0$, there exists $C(p,t,\|b\|_{t,\infty,\beta})>0$ (independent of $u,h,f$ and $g$) such that
\begin{eqnarray}\label{2.7}
\|u\|_{t,\infty,1+\gamma-,p}
\leq C(p,t,\|b\|_{t,\infty,\beta})\[\|h\|_{t,\infty,\beta,p}+ \|f\|_{t,\infty,\alpha} +\|g\|_{t,\infty,p+,E,\alpha}\],
\end{eqnarray}
where
$$
\cC^{1+\gamma-}_b(\mR^d)=\lim_{\varepsilon\rightarrow 0+} \cC^{1+\gamma-\epsilon}_b(\mR^d) =\cap_{0<r<\gamma}\cC^{1+r}_b(\mR^d), \ \ L^{p+}(E,\nu)=\lim_{\epsilon\rightarrow 0+}L^{p+\epsilon}(E,\nu).
$$
\end{theorem}

\begin{remark} \label{rem2.3} (i) Let $k$ be a measurable function on $(E,\nu)$ and $k$ vanish near $0$. For any $1\leq r_1\leq r_2$, if $k\in L^{r_2}(E,\nu)$, then $k\in L^{r_1}(E,\nu)$ and $\|k\|_{r_1,E}\leq C\|k\|_{r_2,E}$. Noticing that $g(t,x,\cdot)$ vanishes near $0$, $g(t,x,\cdot)\in L^{p+}(E,\nu)$, we have $g(t,x,\cdot)\in \cup_{r>p}L^r(E,\nu)$, which implies that there is a positive real number $\epsilon>0$, such that $g(t,x,\cdot)\in L^{p+\epsilon}(E,\nu)$. Therefore, (\ref{2.7}) can be understood as: for every $\epsilon_1>0$, which is sufficiently small, there is a small enough positive real number $\epsilon_2$ $(\epsilon_2\leq\epsilon$), and for every $t>0$, there exists $C(p,t,\|b\|_{t,\infty,\beta})>0$ (independent of $u,h,f$ and $g$) such that
\begin{eqnarray}\label{2.8}
\|u\|_{t,\infty,1+\gamma-\epsilon_1,p}
\leq C(p,t,\|b\|_{t,\infty,\beta})\[\|h\|_{t,\infty,\beta,p}+ \|f\|_{t,\infty,\alpha} +\|g\|_{t,\infty,p+\epsilon_2,E,\alpha}\].
\end{eqnarray}

(ii) When $b=0$,  from the proof, one also asserts that: for every $p\geq 2$ and $g\in L^\infty_{loc}([0,\infty);L^p(E,\nu;\cC^{\alpha}_{b}(\mR^d)))$ with  $g(t,x,\cdot)$ vanishes near $0$, there is a unique mild solution $u$ to (\ref{1.1}). Moreover, $u\in L^\infty_{loc}([0,\infty);\cC^{1+\gamma}_b(\mR^d;L^p(\Omega)))$ and for every $t>0$, there exists $C>
0$, such that
\begin{eqnarray*}
\|u\|_{t,\infty,1+\gamma,p}\leq C(p,t)\[\|h\|_{t,\infty,\beta,p}+ \|f\|_{t,\infty,\alpha} +\|g\|_{t,\infty,p,E,\alpha}\].
\end{eqnarray*}
\end{remark}

\section{Useful lemmas}\label{sec3}
\setcounter{equation}{0}

We now present several lemmas needed for the proof of the main theorem.

\begin{lemma}\label{lem3.1}(Minkowski inequality \cite{Ste})
Assume that $(S_1, \cF_1,\mu_1)$ and $(S_2, \cF_2,\mu_2)$ are two measure spaces and that $G: S_1 \times S_2 \rightarrow \mR$ is measurable. For given real numbers $1\leq p_1\leq p_2$, we also assume that $G\in L^{p_1}(S_1;L^{p_2}(S_2))$.  Then $G\in L^{p_2}(S_2;L^{p_1}(S_1))$ and
\begin{eqnarray}\label{3.1}
\[\int\limits_{S_2}\Big(\int\limits_{S_1}|G(x,y)|^{p_1}\mu_1(dx)\Big)
^{\frac{p_2}{p_1}}\mu_2(dy)\]^{\frac{1}{p_2}}\leq \[\int\limits_{S_1}\Big(\int\limits_{S_2}|G(x,y)|^{p_2}\mu_2(dy)
\Big)^{\frac{p_1}{p_2}}\mu_1(dx)\]^{\frac{1}{p_1}}.
\end{eqnarray}
\end{lemma}

The next lemmas will paly  important roles in estimating stochastic integrals.

\begin{lemma}\label{lem3.2}(Burkholder's inequality \cite[Theorem 4.4.21]{App}) Let $F$ be an $\{\cF_t\}_{t\geq 0}$ adapted stochastic process. Suppose that $\{M_t\}_{t\geq 0}$ is a Brownian type integral of the form
$$
M_t=\int\limits_0^tF(r)dW_r,
$$
for which $F\in L^p(\Omega;L^2_{loc}([0,\infty)))$. Then for any $p\geq 2$, there exists a positive constant $C(p)>0$ such that, for each $t\geq0$,
\begin{eqnarray*}
\mE[|M_t|^{p}] \leq
C(p)\mE\[\int\limits^t_0|F(r)|^2dr\]^{\frac{p}{2}}.
\end{eqnarray*}
\end{lemma}

\begin{corollary}\label{cor3.1} Let $F$ be a $\cB(\mR_+)\times \cB(\mR_+)\times \cB(\mR^d)$-measurable function. Suppose that $\{M_t(x)\}_{t\geq 0}$ is a Brownian type integral of the form
$$
M_t(x)=\int\limits_0^tF(t,r,x)dW_r,
$$
for which
\begin{eqnarray}\label{3.2}
\int\limits^t_0|F(t,r,x)|^2dr<\infty, \quad \mbox{for almost everywhere} \ \ x\in \mR^d.
\end{eqnarray}
Then for any $p\geq 2$, there exists a positive constant $C(p)>0$, which is independent of $x$, such that for each $t\geq0$,
\begin{eqnarray}\label{3.3}
\mE[|M_t(x)|^{p}] \leq
C(p)\[\int\limits^t_0|F(t,r,x)|^2dr\]^{\frac{p}{2}}.
\end{eqnarray}
\end{corollary}
\textbf{Proof.} First, we assume that $F$ has the following form:
\begin{eqnarray}\label{3.4}
F(t,r,x)=\sum_{j=1}^mF_j(t,x)1_{(t_{j-1},t_j]}(r),
\end{eqnarray}
where $m\in\mN$, $F_j$ are $(\mR_+\times\mR^d;\cB(\mR_+)\times \cB(\mR^d))$-measurable, and $0=t_0<t_1<t_2<_{\cdots}<t_m=t$.

Using Lemma \ref{lem3.2} for $p=2$,  we obtain
\begin{eqnarray}\label{3.5}
\mE|M_t(x)|^2=\mE|\sum_{j=1}^m
(W_{t_j}-W_{t_{j-1}})F_j(t,x)|^2
=\sum_{j=1}^m|F_j(t,x)|^2(t_j-t_{j-1})
=\int\limits_0^t|F(t,r,x)|^2dr.
\end{eqnarray}
For $p=4$,  we also  have
\begin{eqnarray}\label{3.6}
&&\mE|M_t(x)|^4\cr&=&\mE|\sum_{j=1}^m
(W_{t_j}-W_{t_{j-1}})F_j(t,x)|^4
\cr&=&\sum_{j=1}^m\mE|
W_{t_j}-W_{t_{j-1}}|^4|F_j(t,x)|^4+6\sum_{i\neq j}\mE|
W_{t_i}-W_{t_{i-1}}|^2\mE|
W_{t_j}-W_{t_{j-1}}|^2|F_i(t,x)|^2|F_j(t,x)|^2
\cr&=&6\[\sum_{j=1}^m|t_j-t_{j-1}|^2|F_j(t,x)|^4+\sum_{i\neq j}(
t_i-t_{i-1})(t_j-t_{j-1})|F_i(t,x)|^2|F_j(t,x)|^2\]
\cr&\leq& 6\[\sum_{j=1}^m|t_j-t_{j-1}||F_j(t,x)|^2\]^2
\cr&=&6\[\int\limits_0^t|F(t,r,x)|^2dr\]^2.
\end{eqnarray}
On the other hand, we have $L^p$-interpolating formulation
\begin{eqnarray}\label{3.7}
F\in L^{p_1}\cap L^{p_3} \  \Longrightarrow \ \|F\|_{L^{p_2}} \leq
\|F\|_{L^{p_1}}^{\frac{(p_3-p_2)p_1}{(p_3-p_1)p_2}} \|F\|_{L^{p_3}}^{\frac{(p_2-p_1)p_3}{(p_3-p_1)p_2}}, \ \  \forall \ \ p_1 \leq p_2
\leq p_3.
\end{eqnarray}
Combining (\ref{3.5}), (\ref{3.6}) and (\ref{3.7})  for   $p\in (2,4)$,  we conclude that there exists $C(p)>0$, which in independent of $m$, such that
\begin{eqnarray}\label{3.8}
\mE|M_t(x)|^p\leq C(p)\[\int\limits_0^t|F(t,r,x)|^2dr\]^{\frac{p}{2}}.
\end{eqnarray}
Observing that the functions  which   satisfy the condition (\ref{3.2}) can be approximated by the step functions of the form (\ref{3.4}), and for $p\in [2,4]$, (\ref{3.2}) holds for step functions,  we thus complete the proof for $p\in [2,4]$.

Analogously, for every even number and every step function of the form (\ref{3.4}), one can prove that (\ref{3.3}) holds. In view of (\ref{3.7}), one derives an inequality of (\ref{3.8}) for every $p>4$. Then by an approximating argument, we complete the proof.

\begin{remark} \label{rem3.1}  When $F(t,r,x)=F(t-r,x)=e^{(t-r)A}f(r,\cdot)(x)$ (with $A$   the generator of a strongly continuous semigroup), we obtain  a Burkholder type inequality for a stochastic convolution. Such an estimate  was   considered by Kotelenez \cite{Kot82} for a square integral martingales with the stochastic convolution taking values in a Hilbert space. The Hilbert space that the Burkholder inequality holds, is then generalized to $2$-uniformly smooth Banach spaces, and so in particular in the Lebesgue spaces $W^{k,q}(\mR^d) \ (2\leq q<\infty)$. This follows from \cite[Proposition 2.4]{Pra}, \cite[Lemma 1.1]{LLP}. Other regularities and related problems for stochastic convolution taking values in $2$-uniformly smooth Banach spaces can be seen \cite{ANP,Da,DaL,HS01,MR,NVW} and the references cited therein. Since $L^\infty(\mR^d)$  is not a $2$-uniformly smooth Banach space, the following inequality in general will be not true:
\begin{eqnarray}\label{3.9}
\mE\[\Big\|\int_0^te^{(t-r)A}f(r,\cdot)dW_r\Big\|_{L^\infty(\mR^d)}^{p}\] \leq
C(p)\[\int\limits^t_0\Big\|\int_0^te^{(t-r)A}f(r,\cdot)\Big\|_{L^\infty(\mR^d)}^2dr\]^{\frac{p}{2}}.
\end{eqnarray}
However, instead of (\ref{3.9}), as a consequence of (\ref{3.3}), we can get
\begin{eqnarray}\label{3.10}
\Big\|\mE\Big|\int_0^te^{(t-r)A}f(r,\cdot)dW_r\Big|^p\Big\|_{L^\infty(\mR^d)} \leq
C(p)\Big\|\[\int\limits^t_0|e^{(t-r)A}f(r,\cdot)|^2dr
\]^{\frac{p}{2}}\Big\|_{L^\infty(\mR^d)}.
\end{eqnarray}
\end{remark}

\begin{lemma}\label{lem3.3} (Kunita's first inequality \cite[Theorem 4.4.23]{App}) Let $E=B_c(0)-\{0\}$ ($0<c\in\mR$). Suppose that $H\in L^p(\Omega;L^2_{loc}([0,\infty);L^2(E,\nu))\cap L^p_{loc}([0,\infty);L^p(E,\nu)))$ is an $\{\cF_t\}_{t\geq 0}$ adapted stochastic process and
\begin{eqnarray*}
I_t=\int\limits_{(0,t]}\int\limits_EH(r,v)\tilde{N}(dr,dv).
\end{eqnarray*}
Then for every $p\geq 2$ and $t\geq 0$, there exists a positive constant $C(p)>0$, such that
\begin{eqnarray*}
\mE[|I_t|^{p}]\leq C(p)\Big\{\mE
\Big[\int\limits^t_0\int\limits_E|H(r,v)|^2\nu(dv)dr\Big]^{\frac{p}{2}}+
\mE\int\limits^t_0\int\limits_E|H(r,v)|^p\nu(dv)dr \Big\}.
\end{eqnarray*}
\end{lemma}

From above Lemma, combining a similar manipulation of Corollary \ref{cor3.1}, one derives that
\begin{corollary}\label{cor3.2} Let $H$ be a $\cB(\mR_+)\times \cB(\mR_+)\times \cB(\mR^d)\times\cB(E)$-measurable function. Suppose that $\{I_t(x)\}_{t\geq 0}$ is a Poisson type integral of the form
\begin{eqnarray}\label{3.11}
I_t(x)=\int\limits_{(0,t]}\int\limits_EH(t,r,x,v)\tilde{N}(dr,dv),
\end{eqnarray}
for which $H(t,r,x,\cdot)$ vanishes near 0 and
\begin{eqnarray*}
\int\limits^t_0\int\limits_E|H(t,r,x,v)|^p\nu(dv)dr<\infty, \ \mbox{for almost everywhere} \ x\in \mR^d.
\end{eqnarray*}
Then for any $p\geq 2$, there exists a positive constant $C(p)>0$, which is independent of $x$, such that for each $t\geq0$,
\begin{eqnarray}\label{3.12}
\mE[|I_t(x)|^{p}]\leq C(p)\int\limits^t_0\int\limits_E|H(t,r,x,v)|^p\nu(dv)dr.
\end{eqnarray}
\end{corollary}
\textbf{Proof.}   Suppose that $H$ has the following form first:
\begin{eqnarray}\label{3.13}
H(t,r,x,v)=\sum_{i=1}^{m_1}\sum_{j=1}^{m_2}H_{i,j}(t,x)1_{(t_{i-1},t_i]}(r)1_{E_j}(v),
\end{eqnarray}
where $m_1,m_2\in\mN$, $H_{i,j}$ are $(\mR_+\times\mR^d;\cB(\mR_+)\times \cB(\mR^d))$-measurable, $0=t_0<t_1<t_2<_{\cdots}<t_{m_1}=t$, $E_j\in\cB(E)$ and $E_{j_1}\cap E_{j_2}=\varnothing$ ($j_1\neq j_2$).

Using Lemma \ref{lem3.3} and the property of independently scattered of Poisson random measure (see Definition \ref{def2.1} (ii)), combining  an analogous argument of Corollary \ref{cor3.1},  we thus have for $p=2$,
\begin{eqnarray}\label{3.14}
\mE[|I_t(x)|^2]=\int\limits^t_0\int\limits_E|H(t,r,x,v)|^2\nu(dv)dr.
\end{eqnarray}
For $p=4$,  we further have
\begin{eqnarray}\label{3.15}
\mE[|I_t(x)|^4]\leq C\Big\{
\Big[\int\limits^t_0\int\limits_E|H(t,r,x,v)|^2\nu(dv)dr\Big]^2+
\int\limits^t_0\int\limits_E|H(t,r,x,v)|^4\nu(dv)dr \Big\}.
\end{eqnarray}

Since $H(t,r,x,\cdot)$ vanishes near $0$, with the help of H\"{o}lder's inequality, from (\ref{3.15})
\begin{eqnarray}\label{3.16}
\mE[|I_t(x)|^4]\leq C\int\limits^t_0\int\limits_E|H(t,r,x,v)|^4\nu(dv)dr.
\end{eqnarray}

From (\ref{3.14}) and (\ref{3.16}), one concludes that: for every  $t>0$, the linear operator given by (\ref{3.11}) is bounded from $L^2([0,t]\times E)$ to $L^2(\Omega)$ and also bounded from $L^4([0,t]\times E)$ to $L^4(\Omega)$. Then the Marcinkiewicz interpolation theorem (\cite[Theorem 2.58]{AF}) is applied, for every $p\in (2,4)$, (\ref{3.12}) holds for step functions, therefore one finishes the proof for $p\in [2,4]$ by an approximating argument.

The remaining part  is the similar as in the proof of    Corollary \ref{cor3.1}. We thus completes the proof.

\begin{remark} \label{rem3.2}  If $H(t,r,x,v)$ is replaced by $U(t,r)h(r-)$ ($U$ is a evolution operator), the estimate was discussed by Kotelenez \cite{Kot84} initially for a square integral martingale with the stochastic convolution taking values in a Hilbert space. The estimate was then strengthened by Ichikawa \cite{Ich}, Hamedani and Zangeneh \cite{HZ}. Some other extensions can also be seen in \cite{BH,DMN,Hau,HS08} and these results are concerned on stochastic evolution taking values in Banach spaces of martingale type $1<p<\infty$. As noticed in  \cite[Remark 2.11]{Hau}, $L^\infty(\mR^d)$ is not a Banach space of martingale type $p$ for any $p>1$, the estimate of (\ref{3.9}) for Poisson random measure in general will be not true. However, as a consequence of (\ref{3.12}), if $H(t,r,x,v)=e^{(t-r)A}h(r,\cdot,v)(x)$ ($A$ is the generator of a strongly continuous semigroup), we can get
\begin{eqnarray}\label{3.17}
&&\Big\|\mE\Big|\int\limits^t_0\!\!\int\limits_Ee^{(t-r)A}h(r,\cdot,v)\tilde{N}(dr,dv)\Big|^p\Big\|_{L^\infty(\mR^d)} \cr&\leq&
C\Big\|\int\limits^t_0\!\!\int\limits_E|e^{(t-r)A}h(r,\cdot,v)|^p\nu(dv)dr
\Big\|_{L^\infty(\mR^d)}.
\end{eqnarray}
The estimate (\ref{3.17}) will play an important role in proving Schauder estimates later.
\end{remark}

\section{Proof of Theorem \ref{the2.1}}\label{sec4}
\setcounter{equation}{0}

\textbf{Proof.}  We divide the proof into three parts: uniqueness, existence and regularity.

\vskip2mm\noindent
\textbf{(Uniqueness).} The stochastic transport-diffusion equation (\ref{1.1}) is linear, to show the uniqueness, it suffices to show that a mild solution with $h=f=g=0$ vanishes identically. When $h=f=g=0$, it becomes a deterministic equation. By virtue of the classical Schauder estimates, it yields that $u=0$, so the mild solution is unique.

To show the existence and regularity, one firstly assumes that $b=0$.

\vskip2mm\noindent
\textbf{(Existence).}  The result follows by using the
explicit formula
\begin{eqnarray}\label{4.1}
u(t,x)=\int\limits_0^t P_{t-r}h(r,\cdot)(x)dr+\int\limits_0^t P_{t-r}f(r,\cdot)(x)dW_r+\int\limits_{(0,t]} \int\limits_EP_{t-r}g(r,\cdot,v)(x)\tilde{N}(dr,dv),
\end{eqnarray}
where $P_t$ is defined by (\ref{2.3}).

By this obvious representation, $u$ meets the properties $(1)$, $(2)$, $(4)$ in Definition \ref{2.3} (for more details, one also refers to \cite{AWZ}). To prove the existence of mild solutions, we need to show that  $u \in L^\infty_{loc}([0,\infty);W^{1,\infty}(\mR^d;L^2(\Omega)))$.

For every $t>0$, from (\ref{4.1}), with the help of (\ref{3.3}) and (\ref{3.12}), one deduces that
\begin{eqnarray*}
\|u\|_{t,\infty,0,2}^2\leq C(t)\[\|h\|_{t,\infty,0,2}^2+ \|f\|_{t,\infty}^2
+\int\limits_0^t\int\limits_E\|g(r,\cdot,z)\|_{\infty}^2\nu(dz)dr\].
\end{eqnarray*}

Now let us verify that $u \in L^\infty_{loc}([0,\infty);W^{1,\infty}(\mR^d;L^p(\Omega)))$. Denote $K(r,x)=\frac{1}{(2\pi r)^{d/2}}e^{-\frac{|x|^2}{2r}}$, if one uses Corollary \ref{cor3.1} and Corollary \ref{cor3.2}, for given $p$, then
\begin{eqnarray}\label{4.2}
&&\mE|u(t,x)|^p\cr&\leq&C(p)\mE\Big|\int\limits_0^t\int\limits_{\mR^d}K(t-r,x-z)h(r,z)dz
dr\Big|^p\cr&&+
C(p)\Big[\int\limits_0^t\Big|\int\limits_{\mR^d}K(t-r,x-z)f(r,z)dz
\Big|^2dr\Big]^{\frac{p}{2}}\cr&&
+C(p)\int\limits_0^t\int\limits_E\Big|\int\limits_{\mR^d}
K(t-r,x-z)g(r,z,v)dz\Big|^p\nu(dv)dr
\cr\cr&\leq& C(p)\[t^p\|h\|_{t,\infty,0,p}^p
+t^{\frac{p}{2}}\|f\|_{t,\infty}^p+t\|g\|_{t,\infty,p,E,0}^p\]
\cr\cr&\leq& C(p,t)\[\|h\|_{t,\infty,0,p}^p
+\|f\|_{t,\infty}^p+\|g\|_{t,\infty,p,E,0}^p\],
\end{eqnarray}
where in the second inequality, we have used Lemma \ref{lem3.1} and in the last inequality, we have used the H\"{o}lder inequality. Moreover, $C(p,t)$ is continuous, non-decreasing in $t$ and $C(p,t)\rightarrow 0$, as $t\rightarrow 0$.

Now let us calculate $|Du|$. For every $1\leq i\leq d$,
\begin{eqnarray*}
&&\partial_{x_i}u(t,x)\cr&=&\int\limits_0^t
dr\int\limits_{\mR^d}\partial_{x_i}K(t-r,x-z)h(r,z)dz+\int\limits_0^t
dW_r\int\limits_{\mR^d}\partial_{x_i}K(t-r,x-z)f(r,z)dz\cr&&+\int\limits_{(0,t]} \int\limits_E\int\limits_{\mR^d}\partial_{x_i}K(t-r,x-z)g(r,z,v)dz\tilde{N}(dr,dv)
\cr&=&\int\limits_0^t
dr\int\limits_{\mR^d}\partial_{x_i}K(t-r,x-z)[h(r,z)-h(r,x)]dz\cr&&+\int\limits_0^t
dW_r\int\limits_{\mR^d}\partial_{x_i}K(t-r,x-z)[f(r,z)-f(r,x)]dz\cr&&+\int\limits_{(0,t]} \int\limits_E\int\limits_{\mR^d}\partial_{x_i}K(t-r,x-z)[g(r,z,v)-g(r,x,v)]dz\tilde{N}(dr,dv).
\end{eqnarray*}
By Corollary \ref{cor3.1} and Corollary \ref{cor3.2},  we therefore conclude that
\begin{eqnarray}\label{4.3}
&&\mE|\partial_{x_i}u(t,x)|^p
\cr&\leq&C(p)\mE\Big|\int\limits_0^t\int\limits_{\mR^d}\partial_{x_i}K(t-r,x-z)
[h(r,z)-h(r,x)]dz
dr\Big|^p\cr&&+C(p)\Big[\int\limits_0^t\Big|\int\limits_{\mR^d}\partial_{x_i}
K(t-r,x-z)[f(r,z)-f(r,x)]dz\Big|^2dr\Big]^{\frac{p}{2}}\cr&&+C(p)\int\limits_0^t\int\limits_E\Big|\int\limits_{\mR^d}
\partial_{x_i}K(t-r,x-z)[g(r,z,v)-g(r,x,v)]dz\Big|^p\nu(dv)dr.
\end{eqnarray}

Observing that $h\in L^\infty_{loc}([0,\infty);\cC^{\beta}_{b}(\mR^d;L^p(\Omega)))$, $f \in L^\infty_{loc}([0,\infty);\cC^{\alpha}_{b}(\mR^d;\mR^d))$ and  $g\in L^\infty_{loc}([0,\infty);L^{p+}(E,\nu;\cC^{\alpha}_{b}(\mR^d)))$, from (\ref{4.3}), we obtain  that
\begin{eqnarray}\label{4.4}
&&\mE|\partial_{x_i}u(t,x)|^p
\cr &\leq&C(p)[h]_{t,\infty,\beta,p}^p \Big|\int\limits_0^t
r^{\frac{\beta-1}{2}}dr \int\limits_{\mR^d}e^{-\frac{|z|^2}{2}}
|z|^\beta dz \Big|^{p}+C(p)[f]_{t,\infty,\alpha}^p \Big[\int\limits_0^t
r^{\alpha-1}dr\Big]^{\frac{p}{2}} \Big[ \int\limits_{\mR^d}e^{-\frac{|z|^2}{2}}
|z|^\alpha dz \Big]^{p} \cr&&+C(p)[g]_{t,\infty,p,E,\alpha}^p
\int\limits_0^tr^{\frac{(\alpha-1)p}{2}}dr\Big[ \int\limits_{\mR^d}e^{-\frac{|z|^2}{2}}
|z|^\alpha dz \Big]^{p}.
\end{eqnarray}
Because of $\alpha,\beta,\gamma>0$, the first two terms in the right hand side of (\ref{4.4}) are finite. Moreover, by (\ref{2.4}), then $\alpha+2/p>1$, so the last term in the right hand side of (\ref{4.4}) is finite as well. $g(t,x,\cdot)$ vanishes near $0$, therefore, for every $t>0$,
\begin{eqnarray}\label{4.5}
\mE|\partial_{x_i}u(t,x)|^p
&\leq& C(p)\[t^{\frac{(\beta+1)p}{2}}[h]_{t,\infty,\beta,p}^p +
t^{\frac{ p}{2}}[f]_{t,\infty,\alpha}^p
+t^{\frac{(\alpha-1)p+1}{2}} [g]_{t,\infty,p,E,\alpha}^p \]
\cr\cr&\leq&  C(p,t)\[ [h]_{t,\infty,\beta,p}^p +
[f]_{t,\infty,\alpha}^p
+[g]_{t,\infty,p,E,\alpha}^p\],
\end{eqnarray}
where $C(p,t)$ is continuous, non-decreasing in $t$ and $C(p,t)\rightarrow 0$, as $t\rightarrow 0$.

\vskip2mm\noindent
\textbf{(Regularity).}  It remains to show the H\"{o}lder estimate for $Du$. We will demonstrate that: $Du\in L^\infty_{loc}([0,\infty);\cC^{\gamma-}_b(\mR^d;L^p(\Omega)))$ and (\ref{2.7}) holds. Observing that $g(t,x,\cdot)$ vanishes near $0$ and $g\in L^\infty_{loc}([0,\infty);L^{p+}(E,\nu;\cC^{\alpha}_{b}(\mR^d)))$, according to Remark \ref{rem2.3} (i), we should prove that for every $\epsilon_1>0$, which is sufficiently small, there is a small enough positive real number $\epsilon_2(\epsilon_1) \ (\epsilon_2\leq\epsilon$), and for every $t>0$, there exists $C(p,t,\|b\|_{t,\infty,\beta})>0$ (independent of $u,h,f$ and $g$) such that
\begin{eqnarray}\label{4.6}
\|Du\|_{t,\infty,\gamma-\epsilon_1,p}
\leq C(p,t,\|b\|_{t,\infty,\beta})\[\|h\|_{t,\infty,\beta,p}+ \|f\|_{t,\infty,\alpha} +\|g\|_{t,\infty,p+\epsilon_2,E,\alpha}\].
\end{eqnarray}

For every $x,y\in\mR^d$ and $1\leq i\leq d$,
\begin{eqnarray*}
&&\partial_{x_i}u(t,x)-\partial_{y_i}u(t,y)\cr
&=&\int\limits_0^t
dr\int\limits_{|x-z|\leq 2|x-y|}\partial_{x_i}K(t-r,x-z)[h(r,z)-h(r,x)]dz
\cr&&-\int\limits_0^t
dr\int\limits_{|x-z|\leq 2|x-y|}\partial_{y_i}K(t-r,y-z)[h(r,z)-h(r,y)]dz
\cr&&+\int\limits_0^t
dr\int\limits_{|x-z|> 2|x-y|}\partial_{y_i}K(t-r,y-z)[h(r,y)-h(r,x)]dz
\cr
&&+\int\limits_0^t
dr\int\limits_{|x-z|> 2|x-y|}[\partial_{x_i}K(t-r,x-z)-\partial_{y_i}K(t-r,y-z)][h(r,z)-h(r,x)]dz
\cr&&+\int\limits_0^t
dW_r\int\limits_{|x-z|\leq 2|x-y|}\partial_{x_i}K(t-r,x-z)[f(r,z)-f(r,x)]dz
\cr&&-\int\limits_0^t
dW_r\int\limits_{|x-z|\leq 2|x-y|}\partial_{y_i}K(t-r,y-z)[f(r,z)-f(r,y)]dz
\cr&&+\int\limits_0^t
dW_r\int\limits_{|x-z|> 2|x-y|}\partial_{y_i}K(t-r,y-z)[f(r,y)-f(r,x)]dz
\cr
&&+\int\limits_0^t
dW_r\int\limits_{|x-z|> 2|x-y|}[\partial_{x_i}K(t-r,x-z)-\partial_{y_i}K(t-r,y-z)][f(r,z)-f(r,x)]dz
\cr&&+\int\limits_{(0,t]} \int\limits_E\int\limits_{|x-z|\leq 2|x-y|}\partial_{x_i}K(t-r,x-z)[g(r,z,v)-g(r,x,v)]dz\tilde{N}(dr,dv)
\cr&&-\int\limits_{(0,t]} \int\limits_E\int\limits_{|x-z|\leq 2|x-y|}\partial_{y_i}K(t-r,y-z)[g(r,z,v)-g(r,y,v)]dz\tilde{N}(dr,dv)
\cr&&+\int\limits_{(0,t]}\int\limits_E\int\limits_{|x-z|> 2|x-y|}\partial_{y_i}K(t-r,y-z)[g(r,y,v)-g(r,x,v)]dz\tilde{N}(dr,dv)\cr&&+
\int\limits_{(0,t]}\!\!\int\limits_E\!\!\int\limits_{|x-z|> 2|x-y|}\!\![\partial_{x_i}K(t-r,x-z)\!-\! \partial_{y_i}K(t-r,y-z)][g(r,z,v)\!-\! g(r,x,v)]dz\tilde{N}(dr,dv)
\cr\cr&=&:I_1(t)\!+\!I_2(t)\!+\!I_3(t)\!+\!I_4(t)\!+\!I_5(t)\!+\!I_6(t)\!+\!I_7(t)\!+\!
I_8(t)\!+\!I_9(t)\!+\!I_{10}(t)\!+\!I_{11}(t)\!+\!
I_{12}(t).
\end{eqnarray*}

Let us estimate $I_1-I_{12}$. To start with, we manipulate the terms $I_1-I_4$.
For convenience of calculations, we set $p_1=2p/(\alpha p-\beta p-\epsilon_1p+2)$, then $1<p_1<p$ and $\beta-1+2/p_1=\gamma-\epsilon_1$.

With the aid of condition (\ref{2.6}) and Lemma \ref{lem3.1}
\begin{eqnarray}\label{4.7}
\mE|I_1(t)|^p&\leq&C(p)\mE\Big|\int\limits_0^t\int\limits_{|x-z|\leq 2|x-y|}|\partial_{x_i}K(t-r,x-z)||h(r,z)-h(r,x)|dzdr\Big|^p
\cr&\leq& C(p)[h]_{t,\infty,\beta,p}^p \Big|\int\limits_0^t\int\limits_{|x-z|\leq 2|x-y|} r^{-\frac{d+1}{2}}e^{-\frac{|x-z|^2}{2r}}|x-z|^\beta dzdr\Big |^p.
\end{eqnarray}
By utilizing the H\"{o}lder inequality and (\ref{3.1}), from (\ref{4.7}), one arrives at
\begin{eqnarray}\label{4.8}
\mE|I_1(t)|^p&\leq& C(p)t^{\frac{(p_1-1)p}{p_1}}[h]_{t,\infty,\beta,p}^p \Big[\int\limits_0^t\Big|\int\limits_{|x-z|\leq 2|x-y|} r^{-\frac{d+1}{2}}e^{-\frac{|x-z|^2}{2r}}|x-z|^\beta dz\Big|^{p_1}dr\Big ]^{\frac{p}{p_1}}
\cr&\leq&C(p,t)[h]_{t,\infty,\beta,p}^p \Big[\int\limits_{|x-z|\leq 2|x-y|} \Big|\int\limits_0^t r^{-\frac{(d+1)p_1}{2}}e^{-\frac{p_1|x-z|^2}{2r}} dr\Big|^{\frac{1}{p_1}} |x-z|^\beta dz\Big ]^p
\cr&\leq&C(p,t)[h]_{t,\infty,\beta,p}^p \Big[\int\limits_{|x-z|\leq 2|x-y|} \Big|\int\limits_0^\infty r^{\frac{(d+1)p_1}{2}-2}e^{-\frac{p_1r}{2}} dr\Big|^{\frac{1}{p_1}} |x-z|^{\beta-d-1+\frac{2}{p_1}} dz\Big ]^p
\cr\cr&\leq&C(p,t)[h]_{t,\infty,\beta,p}^p |x-y|^{(\beta-1+\frac{2}{p_1})p}
\cr\cr&=&C(p,t)[h]_{t,\infty,\beta,p}^p|x-y|^{(\gamma-\epsilon_1)p}.
\end{eqnarray}

Analogously,
\begin{eqnarray}\label{4.9}
\mE|I_2(t)|^p&\leq& C(p,t)[h]_{t,\infty,\beta,p}^p |x-y|^{(\gamma-\epsilon_1) p}.
\end{eqnarray}

For $I_3$, we employ Gauss-Green's formula primarily to gain
\begin{eqnarray}\label{4.10}
I_3(t)=\int\limits_0^t
dr\int\limits_{|y-z|=2|x-y|}K(t-r,y-z)n_i[h(r,y)-h(r,x)]dS.
\end{eqnarray}
From (\ref{4.10}),  owing to the Minkowski and H\"{o}lder inequalities, one ends up with
\begin{eqnarray}\label{4.11}
\mE|I_3(t)|^p&\leq& \Big[\int\limits_0^t
\int\limits_{|y-z|=2|x-y|}K(t-r,y-z)\|h(r,y)-h(r,x)\|_{L^p(\Omega)}dSdr\Big]^p
\cr&\leq&
C(p)t^{\frac{(p_1-1)p}{p_1}}[h]_{t,\infty,\beta,p}^p|x-y|^{\beta p}\Big[\int\limits_0^t
\Big|\int\limits_{|x-z|= 2|x-y|}K(t-r,y-z) dS\Big|^{p_1}dr\Big]^{\frac{p}{p_1}}
\cr&\leq&C(p,t)[h]_{t,\infty,\beta,p}^p|x-y|^{\beta p}\Big[\int\limits_0^t
\Big|\int\limits_{|x-z|= 2|x-y|}r^{-\frac{d}{2}}e^{-\frac{|y-z|^2}{2r}} dS\Big|^{p_1}dr\Big]^{\frac{p}{p_1}}.
\end{eqnarray}

Minkowski's inequality puts into use, from (\ref{4.11}), we achieve
\begin{eqnarray}\label{4.12}
&&\mE|I_3(t)|^p\cr&\leq&C(p,t)[h]_{t,\infty,\beta,p}^p|x-y|^{\beta p}\Big[\int\limits_{|x-z|=2|x-y|} \Big(\int\limits_0^\infty r^{-\frac{p_1d}{2}}e^{-\frac{p_1|y-z|^2}{2r}}dr\Big)^{\frac{1}{p_1}}dS\Big]^{p}
\cr&\leq&C(p,t)[h]_{t,\infty,\beta,p}^p|x-y|^{\beta p}\Big[\int\limits_{|x-z|=2|x-y|} |y-z|^{-d+\frac{2}{p_1}}dz\Big]^{p}\Big(\int\limits_0^\infty r^{\frac{p_1d}{2}-2}e^{-\frac{p_1r}{2}}dr\Big)^{\frac{p}{p_1}}
\cr\cr&\leq&C(p,t)[h]_{t,\infty,\beta,p}^p|x-y|^{(\gamma-\epsilon_1) p}.
\end{eqnarray}

To calculate $I_4$, we use (\ref{3.1}) first, the H\"{o}lder inequality next, (\ref{3.1}) third again, and then acquire
\begin{eqnarray*}
&&\mE|I_4(t)|^p\cr&\leq&
C(p)t^{\frac{(p_1-1)p}{p_1}}[h]_{t,\infty,\beta,p}^p
\Big[\int\limits_0^t
\Big|\int\limits_{|x-z|> 2|x-y|}|x-z|^\beta|\partial_{x_i}K(t\!-\!r,x\!-\!z)\!-\!\partial_{y_i}K(t\!-\!r,y\!-\!z)|dz
\Big|^{p_1}dr\Big]^{\frac{p}{p_1}}
\cr
&\leq&C(p,t)[h]_{t,\infty,\beta,p}^p \Big[\int\limits_{|x-z|> 2|x-y|}|x-z|^\beta\Big(\int\limits_0^t|\partial_{x_i}
K(r,x-z)-\partial_{y_i}K(r,y-z)|^{p_1}dr
\Big)^{\frac{1}{p_1}}dz\Big]^{p}.
\end{eqnarray*}
Notice that $|x-z|>2|x-y|$. So for every $\xi\in [x,y]$,
$$
\frac{1}{2}|x-z| \leq |\xi-z|\leq 2|x-z|.
$$
By virtue of mean value inequality, we have
\begin{eqnarray}\label{4.13}
&&\mE|I_4(t)|^p
\cr&\leq&C(p,t)[h]_{t,\infty,\beta,p}^p |x-y|^p\Big[\int\limits_{|x-z|> 2|x-y|}|x-z|^\beta \Big(\int\limits_0^t r^{-\frac{(d+2)p_1}{2}}e^{-\frac{p_1|x-z|^2}{8r}}dr\Big)^{\frac{1}{p_1}}dz \Big]^{p} \cr&\leq&C(p,t)[h]_{t,\infty,\beta,p}^p |x-y|^p\Big[\int\limits_{|x-z|> 2|x-y|}|x-z|^{\beta-d-2+\frac{2}{p_1}} \Big(\int\limits_0^\infty
r^{\frac{(d+2)p_1}{2}-2}e^{-\frac{p_1r}{8}}dr\Big)^{\frac{1}{p_1}}dz \Big]^{p}\cr\cr&\leq&C(p,t)[h]_{t,\infty,\beta,p}^p|x-y|^{(\gamma-\epsilon_1)p}.
\end{eqnarray}

Let us estimate $I_5-I_8$ and now we set $p_2=2p/(2-\epsilon_1p)$, then $2<p_2<p$ and $\alpha-1+2/p_2=\gamma-\epsilon_1$.

To calculate the term $I_5$, we use (\ref{3.3}) first to derive
\begin{eqnarray}\label{4.14}
\mE|I_5(t)|^p&\leq&C(p)\Big[\int\limits_0^t\Big|\int\limits_{|x-z|\leq 2|x-y|}\partial_{x_i}K(t-r,x-z)[f(r,z)-f(r,x)]dz\Big|^2dr\Big ]^{\frac{p}{2}}
\cr&\leq& C(p)[f]_{t,\infty,\alpha}^p\Big[\int\limits_0^t\Big |\int\limits_{|x-z|\leq 2|x-y|} r^{-\frac{d+1}{2}}e^{-\frac{|x-z|^2}{2r}}|x-z|^\alpha dz\Big |^2dr\Big ]^{\frac{p}{2}}.
\end{eqnarray}
Then the H\"{o}ler inequality applies, for every $t>0$,  we obtain  from (\ref{4.14}) that
\begin{eqnarray*}
\mE|I_5(t)|^p&\leq&C(p)t^{\frac{(p_2-2)p}{2p_2}}[f]_{t,\infty,\alpha}^p \Big[\int\limits_0^t\Big|\int\limits_{|x-z|\leq 2|x-y|} r^{-\frac{d+1}{2}}e^{-\frac{|x-z|^2}{2r}}|x-z|^\alpha dz\Big|^{p_2}dr\Big ]^{\frac{p}{p_2}}.
\end{eqnarray*}

Since $p_2>2$, with the help of the Minkowski inequality, it yields that
\begin{eqnarray}\label{4.15}
&&\mE|I_5(t)|^p\cr&\leq&C(p,t)[f]_{t,\infty,\alpha}^p\Big[\int\limits_{|x-z|\leq 2|x-y|} \Big(\int\limits_0^t r^{-\frac{(d+1)p_2}{2}}e^{-\frac{p_2|x-z|^2}{2r}}dr\Big)^{\frac{1}{p_2}} |x-z|^\alpha dz\Big ]^{p}
\cr
&\leq&C(p,t)[f]_{t,\infty,\alpha}^p\Big[\int\limits_{|x-z|\leq 2|x-y|} \Big(\int\limits_0^\infty r^{\frac{(d+1)p_2}{2}-2}e^{-\frac{p_2r}{2}}drdr\Big)^{\frac{1}{p_2}} |x-z|^{\alpha+\frac{2}{p_2}-d-1} dz\Big ]^{p}
\cr\cr&\leq&C(p,t)[f]_{t,\infty,\alpha}^p|x-y|^{(\alpha-1+\frac{2}{p_2})p}.
\cr\cr&=& C(p,t)[f]_{t,\infty,\alpha}^p|x-y|^{(\gamma-\epsilon_1) p}.
\end{eqnarray}

Similar calculations also imply that
\begin{eqnarray}\label{4.16}
\mE|I_6(t)|^p\leq C(p,t)[f]_{t,\infty,\alpha}^p|x-y|^{(\gamma-\epsilon_1) p}.
\end{eqnarray}

For $I_7$, we employ Gauss-Green's formula firstly to gain
\begin{eqnarray}\label{4.17}
I_7(t)=\int\limits_0^t
dW_r\int\limits_{|y-z|=2|x-y|}K(t-r,y-z)n_i[f(r,y)-f(r,x)]dz
\end{eqnarray}
From (\ref{4.17}), by applying  Corollary \ref{cor3.1} and the H\"{o}lder inequality,  we have that
\begin{eqnarray}\label{4.18}
\mE|I_7(t)|^p&\leq& C(p)\Big[\int\limits_0^t
\Big|\int\limits_{|x-z|= 2|x-y|}K(t-r,y-z)n_i[f(r,y)-f(r,x)]dS\Big|^2dr\Big]^{\frac{p}{2}}
\cr&\leq& C(p)t^{\frac{(p_2-2)p}{2p_2}}
\Big [\int\limits_0^t
\Big|\int\limits_{|x-z|= 2|x-y|}K(t-r,y-z)n_i[f(r,y)-f(r,x)]dS\Big|^{p_2}dr\Big ]^{\frac{p}{p_2}}
\cr&\leq&C(p,t)[f]_{t,\infty,\alpha}^p|x-y|^{\alpha p}\Big [\int\limits_0^t
\Big|\int\limits_{|x-z|= 2|x-y|}r^{-\frac{d}{2}}e^{-\frac{|y-z|^2}{2r}} dS\Big|^{p_2}dr\Big ]^{\frac{p}{p_2}},
\end{eqnarray}
which also suggests that
\begin{eqnarray}\label{4.19}
&&\mE|I_7(t)|^p\cr&\leq&C(p,t)[f]_{t,\infty,\alpha}^p|x-y|^{\alpha p}
\Big[\int\limits_{|x-z|=2|x-y|} \Big(\int\limits_0^\infty r^{-\frac{p_2d}{2}}e^{-\frac{p_2|y-z|^2}{2r}}dr\Big)^{\frac{1}{p_2}}dS\Big]^{p}
\cr&\leq&C(p,t)[f]_{t,\infty,\alpha}^p|x-y|^{\alpha p}\Big[\int\limits_{|x-z|=2|x-y|} |y-z|^{-d+\frac{2}{p_2}}dz\Big]^{p}\Big(\int\limits_0^\infty r^{\frac{p_2d}{2}-2}e^{-\frac{p_2r}{2}}dr\Big)^{\frac{p}{p_2}}
\cr&\leq&C(p,t)[f]_{t,\infty,\alpha}^p|x-y|^{(\alpha+\frac{2}{p_2}-1)p}\Big(\int\limits_0^\infty r^{\frac{p_2d}{2}-2}e^{-\frac{p_2r}{2}}dr\Big)^{\frac{p}{p_2}}.
\cr\cr&\leq&C(p,t)[f]_{t,\infty,\alpha}^p|x-y|^{(\gamma-\epsilon_1)p},
\end{eqnarray}
by using Minkowski's inequality, and in the last inequality, one has used $p_2>2$.

We estimate $I_8$ in terms of Corollary \ref{cor3.1} first, the H\"{o}lder inequality second, the Minkowski inequality third, and then acquire
\begin{eqnarray*}
&&\mE|I_8(t)|^p\cr&\leq&C(p)\Big[\int\limits_0^t
\Big|\int\limits_{|x-z|> 2|x-y|}[\partial_{x_i}K(t-r,x-z)-\partial_{y_i}K(t-r,y-z)][f(r,z)-f(r,x)]dz
\Big|^2dr\Big]^{\frac{p}{2}}
\cr
&\leq&
C(p)t^{\frac{(p_2-2)p}{2p_2}}[f]_{t,\infty,\alpha}^p\Big[\int\limits_0^t
\Big|\int\limits_{|x-z|> 2|x-y|}|x-z|^\alpha|\partial_{x_i}K(t-r,x-z)-\partial_{y_i}K(t-r,y-z)|dz
\Big|^{p_2}dr\Big]^{\frac{p}{p_2}}
\cr
&\leq&C(p,t)[f]_{t,\infty,\alpha}^p \Big[\int\limits_{|x-z|> 2|x-y|}|x-z|^\alpha \Big(\int\limits_0^t|\partial_{x_i}K(r,x-z)-\partial_{y_i}K(r,y-z)|^{p_2}dr
\Big)^{\frac{1}{p_2}}dz\Big]^{p}.
\end{eqnarray*}
Notice that $|x-z|>2|x-y|$, so for every $\xi\in [x,y]$,
$$
\frac{1}{2}|x-z| \leq |\xi-z|\leq 2|x-z|.
$$
In light of mean value inequality, we thus have
\begin{eqnarray}\label{4.20}
&&\mE|I_8(t)|^p
\cr&\leq&C(p,t)[f]_{t,\infty,\alpha}^p |x-y|^p\Big[\int\limits_{|x-z|> 2|x-y|}|x-z|^\alpha \Big(\int\limits_0^t r^{-\frac{(d+2)p_2}{2}}e^{-\frac{p_2|x-z|^2}{8r}}dr\Big)^{\frac{1}{p_2}}dz \Big]^{p}  \cr&\leq&C(p,t)[f]_{t,\infty,\alpha}^p |x-y|^p\Big[\int\limits_{|x-z|> 2|x-y|}|x-z|^{\alpha-d-2+\frac{2}{p_2}} \Big(\int\limits_0^\infty
r^{\frac{(d+2)p_2}{2}-2}e^{-\frac{p_2r}{8}}dr\Big)^{\frac{1}{p_2}}dz \Big]^{p}\cr&\leq&C(p,t)[f]_{t,\infty,\alpha}^p|x-y|^{(\gamma-\epsilon_1) p}.
\end{eqnarray}

Now let us calculate $I_9-I_{12}$. Analogue manipulations for $I_5-I_8$ is applied here, firstly, according to Corollary \ref{cor3.2}, it yields that
\begin{eqnarray}\label{4.21}
&&\mE|I_9(t)|^p
\cr&\leq& C(p)\int\limits_0^t\int\limits_E\Big|\int\limits_{|x-z|\leq 2|x-y|}\partial_{x_i}K(t-r,x-z)[g(r,z,v)-g(r,x,v)]dz\Big|^p\nu(dv)dr
\cr&\leq& C(p)t^\frac{\epsilon_2}{p+\epsilon_2}\Big[\int\limits_0^t\int\limits_E\Big|\int\limits_{|x-z|\leq 2|x-y|}\partial_{x_i}K(t-r,x-z)[g(r,z,v)-g(r,x,v)]dz\Big|^{p+\epsilon_2}\nu(dv)dr
\Big]^{\frac{p}{p+\epsilon_2}}
\cr&\leq&C(p,t)\Big[\int\limits_0^t\int\limits_E\Big|\int\limits_{|x-z|\leq 2|x-y|}|\partial_{x_i}K(t-r,x-z)||x-z|^\alpha dz\Big|^{p+\epsilon_2}[g(r,\cdot,v)]_\alpha^{p+\epsilon_2}\nu(dv)dr
\Big]^{\frac{p}{p+\epsilon_2}}\cr&\leq&C(p,t)[g]_{t,\infty,p+\epsilon_2,E,\alpha}^p
\Big[\int\limits_0^t\Big|\int\limits_{|x-z|\leq 2|x-y|}|\partial_{x_i}K(t-r,x-z)||x-z|^\alpha dz\Big|^{p+\epsilon_2} dr\Big]^{\frac{p}{p+\epsilon_2}}.
\end{eqnarray}

The calculations from (\ref{4.14}) to (\ref{4.15}) applying here again   lead to
\begin{eqnarray}\label{4.22}
\mE|I_9(t)|^p\leq C(p,t)[g]_{t,\infty,p+\epsilon_2,E,\alpha}^p|x-y|^{(\alpha+
\frac{2}{p+\varepsilon_2}-1)p}.
\end{eqnarray}
By setting $\epsilon_2=\epsilon_1p^2/(2-\epsilon_1p)$, then as $\epsilon_1\rightarrow 0$, $\epsilon_2\rightarrow 0$, so if $\epsilon_1$ is small enough, $\epsilon_2<\epsilon$. Moreover, $\alpha+
2/(p+\varepsilon_2)-1=\gamma-\epsilon_1$, thus one concludes that
\begin{eqnarray}\label{4.23}
\mE|I_9(t)|^p\leq C(p,t)[g]_{t,\infty,p+\epsilon_2,E,\alpha}^p|x-y|^{(\gamma-\epsilon_1)p}.
\end{eqnarray}

By the previous argument, we  conclude that
\begin{eqnarray}\label{4.24}
\mE|I_{10}(t)|^p   \vee  \mE|I_{11}(t)|^p \vee \mE|I_{12}(t)|^p \leq C(p)t^\frac{\epsilon_2}{p+\epsilon_2} [g]_{t,\infty,p+\epsilon_2,E,\alpha}^p|x-y|^{(\gamma-\epsilon_1) p}.
\end{eqnarray}

Combining (\ref{4.8})-(\ref{4.9}), (\ref{4.12})-(\ref{4.13}), (\ref{4.15})-(\ref{4.16}), (\ref{4.19})-(\ref{4.20}) and (\ref{4.23})-(\ref{4.24}), we obtain the estimate (\ref{4.6}). According to Remark \ref{rem2.3} (i) (\ref{2.8}) and (\ref{4.2}), (\ref{4.5}), it also implies that
\begin{eqnarray}\label{4.25}
\|u\|_{t,\infty,1+\gamma-,p}\leq C(p,t)\[\|h\|_{t,\infty,\beta,p}+\|f\|_{t,\infty,\alpha}+\|g\|_{t,\infty,p+,E,\alpha}\],
\end{eqnarray}
where $C(p,t)$ is continuous, non-decreasing in $t$ and as $t\rightarrow 0$, $C(p,t)\rightarrow 0$.

From this, we complete the proof for $b=0$. For general $b$, we use the continuity method. First, consider the following family of equations
\begin{eqnarray}\label{4.26}
&&du(t,x)+\theta b(t,x)\cdot\nabla u(t,x)dt-\frac{1}{2}\Delta u(t,x)dt
\cr&&=h(t,x)dt+f(t,x)dW_t+\int\limits_Eg(t,x,v)\tilde{N}(dt,dv), \ \ t>0, \ x\in \mR^d,
\end{eqnarray}
for $\theta\in[0,1]$.  We call a $\theta \in[0,1]$ 'good' if for any \begin{eqnarray}\label{4.27}
\left\{
\begin{array}{ll}
f \in L^\infty_{loc}([0,\infty);\cC^{\alpha}_{b}(\mR^d)), \ h\in L^\infty_{loc}([0,\infty);\cC^{\beta}_{b}(\mR^d;L^p(\Omega))), \\
g\in L^\infty_{loc}([0,\infty);L^{p+}(E,\nu;\cC^{\alpha}_{b}(\mR^d))), \  g(t,x,\cdot) \ \mbox{vanishes near} \ 0,
\end{array}
\right.
\end{eqnarray}
there exists a unique mild solution $u$ to (\ref{4.26}), such that (\ref{2.7}) holds. Moreover if $u\in L^\infty_{loc}([0,\infty);\cC^{1+\gamma-}_b(\mR^d;L^p(\Omega)))$ and (\ref{4.26}) holds in the sense of Definition \ref{2.3} with $b$ in class of $L^\infty_{loc}([0,\infty);\cC^{\beta}_{b}(\mR^d;\mR^d))$ and $f,h,g$ fulfill (\ref{4.27}), by using the estimate (\ref{4.25}), for every given $t>0$, there exists $C(p,t)>0$, which is continuous and non-decreasing in $t$, such that
\begin{eqnarray}\label{4.28}
&&\|u\|_{t,\infty,1+\gamma-,p}
\cr\cr&\leq& C(p,t)\[\|b\cdot\nabla u\|_{t,\infty,\beta,p}+\|h\|_{t,\infty,\beta,p}+\|f\|_{t,\infty,\alpha}
+\|g\|_{t,\infty,p+,E,\alpha}\]
\cr\cr&\leq& C(p,t)\[\|b\|_{t,\infty,\beta} \|Du\|_{t,\infty,\gamma-,p}+\|h\|_{t,\infty,\beta,p}+\|f\|_{t,\infty,\alpha}
+\|g\|_{t,\infty,p+,E,\alpha}\].
\end{eqnarray}
From (\ref{4.28}), in view of fact that $C(p,t)\rightarrow 0$ as $t\rightarrow 0$, for any given $T_0>0$, there is a $T>0$, which is sufficiently small, such that $2C(p,T)<1/\|b\|_{T_0,\infty,\beta}$. Therefore,
\begin{eqnarray}\label{4.29}
\|u\|_{T,\infty,1+\gamma-,p}\leq C(p,T,\|b\|_{T_0,\infty,\beta})
\[\|h\|_{T,\infty,\beta,p}+\|f\|_{T,\infty,\alpha}
+\|g\|_{T,\infty,p+,E,\alpha}\],
\end{eqnarray}
where the constant $C$ is independent of $\theta$. It is clear that $0$ is
a 'good' point.

We now claim that for above $T$, on $[0,T]$ all points of $[0,1]$ are 'good'. To prove the claim we take a 'good' point $\theta_0$ (say $\theta_0=0$) and
rewrite (\ref{4.26}) as
\begin{eqnarray}\label{4.30}
&&du(t,x)+\theta_0 b(t,x)\cdot\nabla u(t,x)dt-\frac{1}{2}\Delta u(t,x)dt
\cr&&=(\theta_0-\theta) b(t,x)\cdot\nabla u(t,x)dt +h(t,x)dt+f(t,x)dW_t+\int\limits_Eg(t,x,v)\tilde{N}(dt,dv).
\end{eqnarray}
For measurable functions $f,g$ and $h$, which satisfy (\ref{4.27}), we define a mapping $S$, which maps $u_1\in L^\infty([0,T];\cC^{1+\gamma-}_{b}(\mR^d;L^p(\Omega)))$  to  the solution $u\in L^\infty([0,T];\cC^{1+\gamma-}_{b}(\mR^d;L^p(\Omega)))$ of the equation
\begin{eqnarray}\label{4.31}
&&du(t,x)+\theta_0 b(t,x)\cdot\nabla u(t,x)dt-\frac{1}{2}\Delta u(t,x)dt
\cr&&=(\theta_0-\theta) b(t,x)\cdot\nabla u_1(t,x)dt +h(t,x)dt+f(t,x)dW_t+\int\limits_Eg(t,x,v)\tilde{N}(dt,dv).
\end{eqnarray}

Observing that owing to our assumptions and the choice of $\theta_0$,
the mapping $S$ is well defined.
Estimate (\ref{4.29}) shows that for any $u_1,\ u_2 \in
L^\infty([0,T];\cC^{1+\gamma-}_{b}(\mR^d;L^p(\Omega)))$
\begin{eqnarray}\label{4.32}
\|Su_1-Su_2\|_{T,\infty,1+\gamma-,p}
\leq C(p,T,\|b\|_{T_0,\infty,\beta})|\theta-\theta_0 |\|u_1-u_2\|_{T,\infty,1+\gamma-,p},
\end{eqnarray}
with $C$ independent of $\theta_0, \theta, u_1$ and $u_2$. It
follows that there is an $\varepsilon> 0$ such that for
$|\theta-\theta_0|\leq \varepsilon$ the mapping $S$ is contractive
in $L^\infty([0,T];\cC^{1+\gamma-}_{b}(\mR^d;L^p(\Omega)))$ and has a
fixed point $u$ which obviously satisfies (\ref{4.26}). Therefore such $\theta$'s are 'good', which certainly proves our claim on time interval $[0,T]$.

Since $u$ is given by (\ref{2.2}),  is  is right continuous in $t$. So $u(T)\in \cC^{1+\gamma-}_{b}(\mR^d;L^p(\Omega))$. We then repeat the proceeding argument   to extend our solution to the time interval $[T,2T]$. Continuing this procedure with
  finitely many steps,  we construct a solution   on the
interval $[0,T_0]$ for every given $T_0>0$. Since $T_0$ is arbitrary, we finish our proof. $\Box$

As seen in the preceding proof, we have a stronger result when the non-Gaussian L\'evy noise is absent ($g=0$).
\begin{corollary} \label{cor4.1} (L\'evy noise is absent: $g=0$)
Consider the stochastic transport-diffusion  equation with Brownian noise
\begin{eqnarray}\label{4.33}
du(t,x)+b(t,x)\cdot\nabla u(t,x)dt-\frac{1}{2}\Delta u(t,x)dt
=h(t,x)dt+f(t,x)dW_t, \ \ t>0, \ x\in \mR^d,
\end{eqnarray}
with the zero initial data.   Assume that $b$ and $h$ satisfy the condition  (\ref{2.6}) with $\beta>0$ and  $p>2$.  Let $f$ satisfy the condition  (\ref{2.5}) with $\alpha>0$. Then there exists a unique mild solution $u$ to (\ref{4.33}) in the space $L^\infty_{loc}([0,\infty);\cC^{1+\alpha-}_b(\mR^d;L^p(\Omega)))$. Moreover for every given $t>0$, there exists $C(p,t,\|b\|_{t,\infty,\beta})>0$ (independent on $u,h$ and $f$) such that
\begin{eqnarray}\label{4.34}
\|u\|_{T,\infty,1+\alpha-,p} \leq C(p,t,\|b\|_{t,\infty,\beta})\[\|h\|_{t,\infty,\beta,p}+ \|f\|_{t,\infty,\alpha}\].
\end{eqnarray}
\end{corollary}
\textbf{Proof.}
If we can prove the case of $b=0$, then by using the continuity method,  we get the conclusion for $b\neq 0$ which fulfills the assumption (\ref{2.6}). Hence, it is sufficient to prove this corollary in the case of  $b=0$.

The existence and uniqueness of mild solution in space $L^\infty_{loc}([0,\infty);W^{1,\infty}(\mR^d;L^p(\Omega)))$ can be deduced from  (\ref{4.2}) and (\ref{4.5}).  Thus it remains to derive  the $\cC^{1+\alpha-}$ estimate.

Using the notations   in   the proof in Theorem \ref{the2.1}, we should prove that for every $\epsilon_1>0$, which is sufficiently small, there is a small enough real positive number $\epsilon_2(\epsilon_1) \ (\epsilon_2\leq\epsilon$), and for every $t>0$, there exists $C(p,t)>0$ (continuous and non-decreasing in $t$, and independent of $u,h,f$) such that
\begin{eqnarray}\label{4.35}
\|Du\|_{t,\infty,\alpha-\epsilon_1,p}
\leq C(p,t)\[\|h\|_{t,\infty,\beta,p}+\|f\|_{t,\infty,\alpha}\].
\end{eqnarray}

For $1< q< p$, we get from  (\ref{4.7}) and (\ref{4.8}),
\begin{eqnarray}\label{4.36}
\mE|I_1(t)|^p\leq C(p)t^{\frac{(q-1)p}{q}}[h]_{t,\infty,\beta,p}^p |x-y|^{(\beta-1+\frac{2}{q})p}.
\end{eqnarray}
Similarly, we arrive at
\begin{eqnarray}\label{4.37}
\mE|I_2(t)|^p \vee \mE|I_3(t)|^p \vee \mE|I_4(t)|^p \leq C(p)t^{\frac{(q-1)p}{q}}[h]_{t,\infty,\beta,p}^p |x-y|^{(\beta-1+\frac{2}{q})p}.
\end{eqnarray}

Take $q=2/(1+\alpha-\beta -\epsilon_1)$. Then $1<q<2<p$ and $\beta-1+2/q=\alpha-\epsilon_1$. From (\ref{4.36}) and (\ref{4.37}), we get
\begin{eqnarray}\label{4.38}
\mE|I_1(t)|^p\vee \mE|I_2(t)|^p \vee \mE|I_3(t)|^p \vee \mE|I_4(t)|^p \leq C(p,t)[h]_{t,\infty,\beta,p}^p|x-y|^{(\alpha-\epsilon_1) p}.
\end{eqnarray}

To calculate $I_5-I_8$, we use (\ref{4.14})-(\ref{4.20}).  For every $2<q_1<p$, we  have
\begin{eqnarray}\label{4.39}
\mE|I_5(t)|^p \vee \mE|I_6(t)|^p \vee \mE|I_7(t)|^p \vee \mE|I_8(t)|^p \leq C(p)t^{\frac{(q_1-2)p}{2q_1}}[f]_{t,\infty,\alpha}^p |x-y|^{(\alpha-1+\frac{2}{q_1})p}.
\end{eqnarray}
  Taking $q_1=2/(1-\epsilon_1)$, then $2<q_1<p$ and
\begin{eqnarray}\label{4.40}
\mE|I_5(t)|^p \vee \mE|I_6(t)|^p \vee \mE|I_7(t)|^p \vee \mE|I_8(t)|^p \leq C(p,t)[f]_{t,\infty,\alpha}^p |x-y|^{(\alpha-\epsilon_1) p}.
\end{eqnarray}
From (\ref{4.38}), (\ref{4.40}) and Remark \ref{rem2.3} (i),  we complete the proof. $\Box$

\section{Concluding remarks}\label{sec5}
\setcounter{equation}{0}

In view of Remark \ref{rem2.3} (ii) and Corollary \ref{cor4.1}, we have seen the following result.  For every $\alpha>0$,  assume that  $f \in L^\infty_{loc}([0,\infty);\cC^{\alpha}_{b}(\mR))$,
and  that  there is a real number $\beta$ (with $0<\beta<\alpha$ )  such that $h\in L^\infty_{loc}([0,\infty);\cC^{\beta}_{b}(\mR;L^2(\Omega)))$. Then the   Cauchy problem
\begin{eqnarray}\label{5.1}
du(t,x)-\frac{1}{2}\Delta u(t,x)dt=h(t,x)dt+f(t,x)dW_t, \ \ t>0, \ x\in \mR, \ \ u|_{t=0}=0,
\end{eqnarray}
has a unique mild solution $u\in L^\infty_{loc}([0,\infty);\cC^{1+\alpha}_b(\mR;L^2(\Omega)))$, which is given by
\begin{eqnarray}\label{5.2}
u(t,x)=\int\limits_0^tdr\int\limits_{\mR}K(t-r,x-z)h(r,z)dz
+\int\limits_0^t
dW_r\int\limits_{\mR}K(t-r,y-z)f(r,z)dz.
\end{eqnarray}
In $L^\infty_{loc}([0,\infty);\cC^{1+\alpha}_b(\mR))$, we   choose a non-negative and time independent function $f$ such that

(i) $f$ in non-decreasing on $\mR$ and $\mbox{supp}f\subset\mR_+$;

(ii) for $x,y\in[0,1]$, $|f(x)-f(y)|\approx |x-y|^\alpha$.

For this function and for  $0<x<1$,  we conclude    by using (\ref{3.5}) that,
\begin{eqnarray}\label{5.3}
&&\mE|\partial_xu(t,x)-\partial_xu(t,0)|^2\cr&=&\int\limits_0^t\Big|\int\limits_{\mR}
\partial_zK(r,z)[f(x-z)-f(-z)]dz\Big|^2dr
\cr&=&\int\limits_0^t\Big|\int\limits_{\mR}
\frac{z}{r}K(r,z)[f(x+z)-f(z)]dz\Big|^2dr
\cr&=&\int\limits_0^t\Big|\int\limits_0^\infty
\frac{z}{r}K(r,z)[f(x+z)-f(z)]dz\Big|^2dr.
\end{eqnarray}

Observe that $f$ is non-negative and non-decreasing. For every $\delta$ (with  $1>\delta>\alpha$), it yields from (\ref{5.3}) that,
\begin{eqnarray*}
&&\sup_{0<x<1}\frac{\|\partial_xu(t,x)-\partial_xu(t,0)\|_2}{x^\delta}
\cr&=&\sup_{0<x<1} \Big[\int\limits_0^t\Big|\int\limits_0^\infty
\frac{z}{r}K(r,z)\frac{f(x+z)-f(z)}{x^\delta}dz\Big|^2dr\Big]^{\frac{1}{2}}
\cr&\geq &\sup_{0<x<\frac{1}{2}} \Big[\int\limits_0^t\Big|\int\limits_0^{\frac{1}{2}}
\frac{z}{r}K(r,z)\frac{f(x+z)-f(z)}{x^\delta}dz\Big|^2dr\Big]^{\frac{1}{2}}
\cr\cr &\geq &C\sup_{0<x<\frac{1}{2}} \Big[\int\limits_0^t\Big|\int\limits_0^{\frac{1}{2}}
\frac{z}{r}K(r,z)\frac{x^\alpha}{x^\delta}dz\Big|^2dr\Big]^{\frac{1}{2}}
\cr\cr&=&\infty.
\end{eqnarray*}
On the other hand, if we chooses $h(t,x,\omega)=h_1(t,x)h_2(\omega)$ with $h_1\in L^\infty_{loc}([0,\infty);\cC^{\beta}_{b}(\mR))$ $h_2\in L^2(\Omega)$,  then the first term in the right hand side of (\ref{5.2}) belongs to $L^\infty_{loc}([0,\infty);\cC^{1+\delta}_b(\mR^d;L^2(\Omega)))$ for every $\delta$ (with  $\alpha<\delta<1$). From this result, we see  that for every $\delta>\alpha$, $u$ is not in class of $L^\infty_{loc}([0,\infty);\cC^{1+\delta}_b(\mR^d;L^2(\Omega)))$, i.e. the H\"{o}lder index $\alpha$ is optimal now.

Analogously, by taking $g(t,x,v)=f(t,x)g_1(v)$, with $f$ satisfying the properties described above and $g_1\in L^2(E,\nu)$,  we claim that:
For every $\alpha>0$, if there is a real number $0<\beta<\alpha$, such that $h\in L^\infty_{loc}([0,\infty);\cC^{\beta}_{b}(\mR;L^2(\Omega)))$, then the following Cauchy problem
\begin{eqnarray}\label{5.4}
du(t,x)-\frac{1}{2}\Delta u(t,x)dt=h(t,x)dt+\int\limits_Eg(t,x,v)\tilde{N}(dt,dv), \ \ t>0, \ x\in \mR, \ \ u|_{t=0}=0,
\end{eqnarray}
has a unique mild solution $u\in L^\infty_{loc}([0,\infty);\cC^{1+\alpha}_b(\mR;L^2(\Omega)))$, which is given by
\begin{eqnarray}\label{5.5}
u(t,x)=\int\limits_0^tdr\int\limits_{\mR}K(t-r,x-z)h(r,z)dz
+\int\limits_0^t \int\limits_E
\int\limits_{\mR}K(t-r,y-z)g(r,z,v)dz\tilde{N}(dt,dv).
\end{eqnarray}
In view of (\ref{3.14}) and (\ref{5.3}), and  using (\ref{3.5}),   we conclude that
for $0<x<1$,
\begin{eqnarray}\label{5.6}
\mE|\partial_xu(t,x)-\partial_xu(t,0)|^2=\int\limits_0^t\Big|\int\limits_0^\infty
\frac{z}{r}K(r,z)[f(x+z)-f(z)]dz\Big|^2dr\int\limits_E|g_1(v)|^2\nu(dv),
\end{eqnarray}
which says that the H\"{o}lder index $\alpha$ is optimal for (\ref{5.4}).

As the stochastic processes $\{W_t\}_{t\geq 0}$ and $\{\tilde{N}_t\}_{t\geq 0}$ are independent, the Wiener-It\^{o} integral in (\ref{5.2}) (interpreted as a stochastic process) and the Wiener-L\'{e}vy integral in (\ref{5.5}) are independent as well. Combining (\ref{5.3}) and (\ref{5.6}), we conclude that when $p=2$, H\"{o}lder index $\alpha$ for the Cauchy problem
\begin{eqnarray*}
du(t,x)\!-\!\frac{1}{2}\Delta u(t,x)dt\!=\!h(t,x)dt\!+\!f(t,x)dW_t\!+\!\int\limits_Eg(t,x,v)\tilde{N}(dt,dv), \ t>0, x\in \mR, \ u|_{t=0}=0,
\end{eqnarray*}
is also optimal.

\vskip4mm\noindent
\textbf{\large{Acknowledgements}}
 \vskip3mm\par
This research was partly supported by the NSF of China grants 11501577, 11301146, 11531006, 11371367  and  11271290.

\end{document}